# Average Exit Time for Volume Preserving Maps

J.D. Meiss

Program in Applied Mathematics

University of Colorado

Boulder, CO 80304

April 5, 1996

**Abstract**

For a volume preserving map, the exit time, averaged over the incoming set of a region, is given by the ratio of the measure of the accessible subset of the region to that of the incoming set. This result is primarily of interest to show two things: first it gives a simple bound on the algebraic decay exponent of the survival probability. Second, it gives a tool which permits the computation of the measure of the accessible set. We use this to find the measure of the bounded orbits for Hénon's quadratic map.

## 1 Introduction

Orbits of an area preserving map are eternally trapped in a region if they are enclosed by an invariant circle, but can leak through destoyed circles (cantori). However, even untrapped orbits feel the influence of these invariant boundaries because they are extremely "sticky."[1-13] In higher dimensional maps, the mechanisms for trapping and escape are much less understood [14]. In this paper we study the time of escape for orbits that begin in a specified region, $A$. We discuss the *exit* time, the time for a point to first exit the set, and the *transit* time, the time for a point to traverse the set. The *exit time distribution*. is the probability distribution of exit times. Our primary goal is to use this distribution to probe the trapped invariant set.

The study of escape time for chaotic maps has long been of interest. As an example, consider Hénon's area preserving, quadratic map, which we write as

$$H: (x,y) \to (y - k + x^2, -x) \quad . \qquad (1)$$

We are interested in the set of bounded orbits. It is possible to show that (for $k>-1$) all bounded orbits are contained in the square $\{ z = (x,y): |x,y| < x_s\}$ where $x_{s,e} = 1 \pm \sqrt{1+k}$ . Here we denote the two fixed points (for $k>-1$) by $z_e = (x_e, -x_e)$ and $z_s = (x_s, -x_s)$. The first is elliptic when $-1<k<3$, and the second is a saddle. Bounded orbits include all periodic orbits and all orbits within invariant circles; the latter exist, for example, in the neighborhood of $z_e$ providing it is elliptic and has a rotation number, ω, which is not $1/3$ or $1/2$ (i.e., $k$ is not $5/4$ or larger than 3). There are also trapped hyperbolic orbits, for example those that are homoclinic to the saddle fixed point.

In his original paper [15], Hénon studied the set of bounded orbits of $H$ (in a different coordinate system than ours) by iterating points along a segment of the symmetry line ($x=-y$) to determine the subintervals that did not reach a fixed large distance from the origin within 100 iterates. He noted that as $k$ varies, the boundaries of the trapped intervals either closely follow the position of an unstable rotational periodic point or else an island around an elliptic periodic point. Channon and Leibowitz [16] studied the escape and trapping from the period 5 island chain in the Hénon map at $k = -0.422$. They identified the exit and entrance lobes of the resonance as the important sets to consider. Studying the $5^{th}$ power of the map, they started 7750 orbits in the outer entrance lobe and computed the survival probability distribution. This was found to decay as $t^{-\alpha}$ with $\alpha = 0.5$ for short times (up to 10 iterates), but deviated from this for moderate times (up to 45 iterates).

Karney [17] also considered the Hénon map, and studied the survival probability for a square $A$ enclosing all bounded orbits. He mostly studied $k= -0.6$ where the most prominent island chain is period 6. Karney's primary object of study was the "trapping time statistic" which is proportional to the exit time distribution for $A$. He found that this distribution decayed as $t^{-\alpha+2}$ with $\alpha$ about 0.25, for times up to $10^4$ iterates, though the slope appears to vary and not settle down to a fixed value for times up to $10^8$. Chirikov and Shepelyansky [18] did similar computations for the standard map for the parameter value when there is a critical golden circle. They obtained $\alpha = 0.34$ up to $10^5$ iterates. The exponent for the decay is apparently not universal [19]. One reason this might be so, according to Murray [20], is that the self-similar limit is not reached for "short" time computations.

Rom-Kedar and Wiggins have emphasized the fact that one can obtain a complete description of transport through a region by considering the future history of only the entering trajectories [5]. Using their notion of "lobe dynamics" for a homoclinic tangle, they obtained an expression for the accessible region in terms of the exit time distribution—we will use a slightly generalized version of this below.

## 2 Definitions

Let $f: M \to M$ be a homeomorphism with an invariant measure $\mu$. For measurable sets $A, B \subset M$, the *crossing time* (or first passage time) for $a \in A$ to $B$ is defined as

$$t_{A \to B}(a) = \min_{n > 0}\left(n : f^n(a) \in B\right). \tag{2}$$

We let $t_{A \to B}(a) = \infty$ if $a$ never reaches $B$. If $B$ is the complement of $A$, then the crossing time is called the (forward) *exit time*

$$t^+(a) = t_{A \to M \setminus A}(a) . \tag{3}$$

Similarly the backward exit time for $a \in A$ is defined as

$$t^-(a) = \min_{n > 0} \left( n : f^{-n}(a) \notin A \right) . \tag{4}$$

If $B = A$, then the crossing time is called the first return time:

$$t_{return}(a) \equiv t_{A \to A}(a) , \text{ for } a \in A. \tag{5}$$

According to our definition $t_{return}(a) \geq 1$, and is equal to 1 whenever the point does not leave on the first iterate. We define the exit set $E \subset A$ as the subset with exit time 1, and the entry set $I \subset A$ as the set with backward exit time 1; equivalently,

$$E = A \setminus f^{-1}(A) , I = A \setminus f(A) . \tag{6}$$

The *turnstile* is the union of $E$ and $I$.

As an example, consider the Hénon map (1). A natural choice for the region $A$ is the "resonance zone", shown in Fig. 1, consisting of the region bounded by the segments of the left-going branches of the stable and unstable manifolds of the saddle fixed point up to their first intersection on the symmetry line ($x=-y$). This region contains all bounded orbits of $H$.

The transit time of a point $a$ is the sum of its forward and backward exit times minus one (to get rid of an annoying term),

$$t_{transit}(a) \equiv t^+(a) + t^-(a) - 1 . \tag{7}$$

The transit time is an orbit invariant: each point along an orbit has the same transit time. The *transit time decomposition* of a set is the partition of a set into subsets with equal transit time. We denote the transit time decomposition of the entry set $I$ by sets $T_j$; this is also the part of $I$ with exit time $j$:

$$T_j = \{ a \in I : t^+(a) = j \} . \tag{8}$$

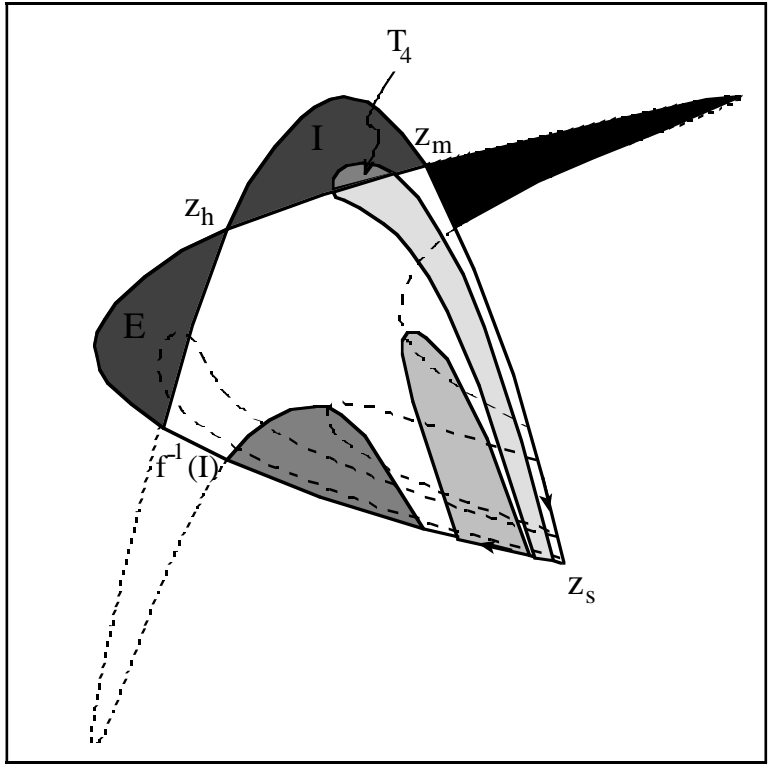

Figure 1. Exit and entry sets for the fixed point resonance of the Hénon map for $k = 0.5$.

For example, we show the first few preimages of the exit set, in sequentially lighter shades, for the Hénon resonance zone in Fig. 1. For this parameter value $T_1 = T_2 = T_3 = \phi$ and $T_4$ is shown. A partial transit time decomposition of the Hénon entry set is shown in Fig. 2.

(color figure attached)

Figure 2. Exit time decomposition of the entry set of the Hénon Trellis at $k=0.5$. Color scale is given at the bottom.

By measure preservation, almost every point that enters $A$ must eventually escape, so

$$\mu(E) = \mu(I) = \sum_{j=1}^{\infty} \mu(T_j) \qquad (9)$$

The accessible set, $A_{acc} \subset A$ is defined to be the set with finite backward exit time; it is the set that can be reached from the outside:

$$A_{acc} = \left\{ a \in A : t^-(a) < \infty \right\} \qquad (10)$$

Of course, the set with finite exit time differs from $A_{acc}$ at most by a set of measure zero, since the set that enters but never exits must have measure zero.

## 3 Average Exit Time

In this note we are interested in these transport times averaged over sets of initial conditions. For a function $g(z)$ we denote

$$\langle g \rangle_S \equiv \frac{1}{\mu(S)} \int_S g(z) \, d\mu$$

Remarkably, there are some simple formulae for average transport times. The following lemma was stated for volume preserving flows (without its elementary proof) in [21] and is implicitly obtained for two dimensional maps in [22].

**Average Exit Time Lemma**: The average exit (transit) time for incoming orbits is

$$\langle t^+ \rangle_I = \langle t_{transit} \rangle_I = \frac{\mu(A_{acc})}{\mu(I)} \tag{11}$$

Proof: Since $T_j \subset I$ are disjoint, and cover almost all of $I$, and the exit (and transit) time of the set $T_j$ is $j$, the average exit time is given by

$$\langle t^+ \rangle_I = \frac{1}{\mu(I)} \sum_{j=1}^{\infty} j \, \mu(T_j) \ ,$$

(assuming the sum exists). Given the transit time decomposition $T_j$, we can compute the accessible subset of $A$ using the sets

$$T_j^i = f^i(T_j) \ ,$$

Note that $T_j^i \subset A$ for $i = 0,...,j-1$, and that these sets are disjoint because they have transit time $j$ and backward exit time $i+1$. Furthermore the union of these sets is, up to measure zero, the entire subset of $A$ that exits. Thus by definition

$$A_{acc} = \bigcup_{j=1}^{\infty} \bigcup_{i=0}^{j-1} T_j^i \quad \Rightarrow \quad \mu(A_{acc}) = \sum_{j=1}^{\infty} \sum_{i=0}^{j-1} \mu(T_j^i) \ .$$

Since measure is preserved $\mu(T_j^i) = \mu(T_j)$, and so

$$\mu(A_{acc}) = \sum_{j=1}^{\infty} j \mu(T_j) \tag{12}$$

This expression was first obtained by Rom Kedar and Wiggins (see e.g., (6.13) in [22]). Finally since $\mu(A_{acc}) \leq \mu(A)$, it is clear that the sum converges. ❏

Since $\mu(A_{acc})$ is finite, a simple consequence of this result is

**Corollary 1:** The measure of the region with transit time t must decay at least as $t^{-2}$.

Furthermore, the lemma implies a well known result of Kac in ergodic theory [23]. In our context this can be generalized to

**Corollary 2 (Kac):** Suppose $\mu(M) = 1$. The average first return time to a region $A \subset M$ is

$$\langle t_{return} \rangle_A = \frac{\mu(M_{acc})}{\mu(A)} , \qquad (13)$$

where $M_{acc}$ is the subset of $M$ that is accessible to orbits beginning in $A$.

Proof: For points that stay in $A$ for at least one step, the first return time is one. The remaining portion of $A$ is its exit set $E$. Now, $f(E)$ is the incoming set of $M\backslash A$, and by our Lemma, the average transit time for these points through $M\backslash A$ is

$$\langle t_{f(E) \to A} \rangle_{f(E)} = \frac{\mu(M_{acc}\backslash A)}{\mu(E)} = \frac{\mu(M_{acc}) - \mu(A)}{\mu(E)} .$$

Since there has been one step already, the first return time to $A$ for these points is one larger. So the average first return time to $A$ is

$$\langle t_{return} \rangle_A = \frac{1}{\mu(A)} \left[ \left( \mu(A) - \mu(E) \right) \times 1 + \mu(E) \times \left( \langle t_{f(E) \to A} \rangle_{f(E)} + 1 \right) \right] .$$

This reduces to the promised result. ❏

We can also easily compute the transit time averaged over $A$—since each image of $T_j$ has the same transit time, the area of $A$ that has transit time $j$ is $j\mu(T_j)$, therefore

**Corollary 3:** The average transit time for points that do escape from $A$ is

$$\langle t_{transit} \rangle_{A_{acc}} = \frac{1}{\mu(A_{acc})} \sum_{j=1}^{\infty} j^2 \mu(T_j) , \qquad (14)$$

providing that this sum converges. The average exit time for $A$ is

$$\langle t^+ \rangle_{A_{acc}} = \frac{1}{\mu(A_{acc})} \sum_{j=0}^{\infty} \mu(T_j) \sum_{k=0}^{j-1} (j-k) = \frac{1}{\mu(A_{acc})} \sum_{j=0}^{\infty} \frac{j(j+1)}{2} \mu(T_j)$$
$$= \frac{1}{2} \left( \langle t_{transit} \rangle_{A_{acc}} + 1 \right) \qquad (15)$$

In contrast to the exit time averaged over *I*, we cannot use disjointness to show that (14) or (15) converges. In fact as we show by a simple example in the next section, these sums need not converge.

## 4 Exit Time Distributions

As was emphasized by Rom-Kedar and Wiggins [5], if one knows the $T_j$, one has most of the information one could want about transport though *A*. Here we recall some definitions of normalized exit and transit time distributions.

The exit time probability distribution for the entry set is the probability that a trajectory in *I* will have exit time *t*:

$$Prob(\, t^+(I) = j) = \frac{\mu(T_j)}{\mu(I)} \,. \tag{16}$$

Similarly the survival probability, is the probability that the exit time will be at least than *k*:

$$Prob(\, t^+(I) \geq k) \ = \frac{1}{\mu(I)} \sum_{j=k}^{\infty} \mu(T_j) \tag{17}$$

Note that $P(t^+(I) \geq 1) = 1$ by eq. (9).

Once we know the $T_j$, distributions for *A* are also known. For example, since $f^j(T_{k+j}) \subset A_{acc}$ has exit time *k* (with backward exit time *j*+1), the subset of $A_{acc}$ that has exit time *k* is given by $\bigcup_{j=0}^{\infty} f^j(T_{k+j})$. Thus the exit time probability distribution for $A_{acc}$ is

$$Prob(\, t^+(A_{acc}) = k) = \frac{1}{\mu(A_{acc})} \sum_{j=k}^{\infty} \mu(T_j) = \frac{1}{<t^+>_I} Prob(t^+(I) \geq k)$$

which is the same as the survival distribution for *I* (17), up to normalization. Similarly the subset of *A* with transit time *j* is $\bigcup_{k=0}^{j} f^k(T_j)$, thus the transit time probability is

$$Prob(\, t_{transit}(A_{acc}) = j) \ = \frac{1}{\mu(A_{acc})} \sum_{k=0}^{t-1} \mu(T_j) = j\frac{\mu(T_j)}{\mu(A_{acc})}$$

Finally the survival time distribution for A is

$$Prob(\, t^+(A) \geq k) = \frac{1}{\mu(A_{acc})} \sum_{j=k}^{\infty} \sum_{m=j}^{\infty} \mu(T_m) = \frac{1}{\mu(A_{acc})} \sum_{j=1}^{\infty} j\mu(T_{k+j})$$

Note that these equations imply that if , e.g.,

$$\mu(T_k) \sim k^{-(2+\alpha)}, \ \text{as } k \to \infty,$$

where $\alpha > 0$ then we have

$$Prob(t^+(I) = k) \sim k^{-(\alpha+2)}$$
$$Prob(t^+(I) \geq k) \sim Prob(t^+(A_{acc}) = k) \sim Prob(t_{transit}(A_{acc}) = k) \sim k^{-(\alpha+1)}$$
$$Prob(t^+(A_{acc}) \geq k) \sim k^{-\alpha}$$

## 5 Examples

Consider the linear, area preserving, hyperbolic map

$$\begin{pmatrix} x' \\ y' \end{pmatrix} = \begin{pmatrix} \lambda & 0 \\ 0 & \lambda^{-1} \end{pmatrix} \begin{pmatrix} x \\ y \end{pmatrix}$$

where $\lambda > 1$. Let A be the unit square $\{(x,y) : 0 \leq x,y \leq 1\}$. Then the entrance set is the rectangle $I = A \setminus f(A) = \{(x,y) : 0 \leq x \leq 1, \lambda^{-1} < y \leq 1\}$. It is easy to see that the transit time decomposition of $I$ is

$$T_j = \{(x,y) : \lambda^{-j} < x \leq \lambda^{1-j}, \lambda^{-1} < y \leq 1\}$$

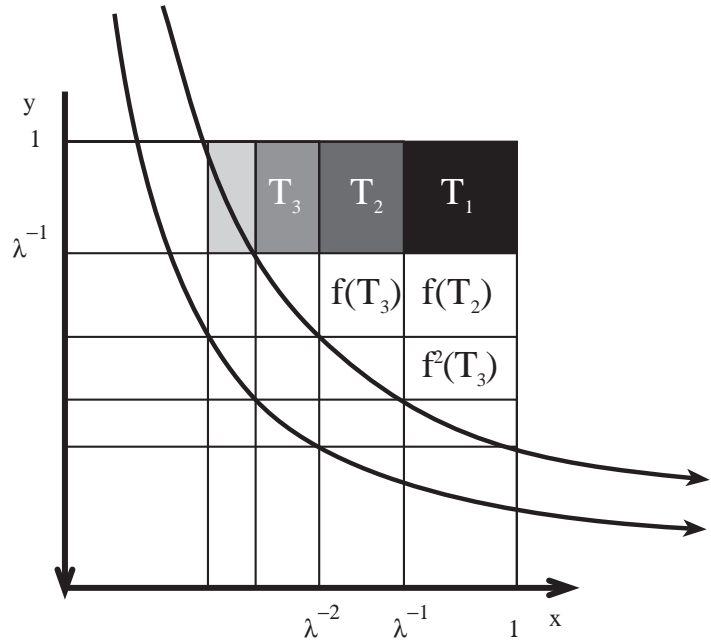

Figure 3. Transit time decomposition of the entry region for the linear hyperbolic map. The transiting regions $T_j$, j=1,2,3 and a few of their iterates are shown.

So the measures of each of these regions are

$$\mu(T_j) = \frac{(\lambda - 1)(1 - \lambda^{-1})}{\lambda^j} \ . \tag{18}$$

These decay exponentially, as one would expect. The calculations needed for Eqs. (9), (12) and (14) are derivatives of simple geometric sums, yielding

$$\sum_{j=1}^{\infty} \mu(T_j) = (1 - \lambda^{-1}) = \mu(I) ,$$

$$\sum_{j=1}^{\infty} j \mu(T_j) = 1 = \mu(A) ,$$

$$\sum_{j=1}^{\infty} j^2 \mu(T_j) = \frac{\lambda + 1}{\lambda - 1} ,$$

as required by the lemma. Thus the average transport times are

$$\langle t^+ \rangle_I = \frac{\lambda}{\lambda - 1}, \quad \langle t_{transit} \rangle_A = \frac{\lambda + 1}{\lambda - 1}, \quad \langle t^+ \rangle_A = \frac{\lambda}{\lambda - 1} . \tag{19}$$

These results are unchanged if we scale the size of $A$ (since the map is linear), or if we replace $A$ by a square centered on the fixed point (since the map is symmetric).

More generally, suppose the map is a diagonal hyperbolic matrix with eigenvalues ($\lambda_1, \lambda_2..., \lambda_d, \mu_{d+1}, \mu_2..., \mu_n$), where the $\lambda_i > 1$ and the $\mu_j < 1$. Then the entry set of the unit hypercube $A$ is the union of rectangular cylinders:

$$I = \bigcup_{k=d+1}^{n} \left\{ (x_1,...x_n) \in A : \mu_k < x_k \leq 1 \right\} .$$

If we define $\Lambda = \prod_{i=1}^{d} \lambda_i$, and $\Pi = \prod_{i=d+1}^{n} \mu_i$ to be the total expansion and contraction, respectively, then a simple calculation gives the measure of the transit regions by

$$\mu(T_j) = \frac{(\Lambda - 1)(1 - \Pi)}{\Lambda^j} , \tag{20}$$

in complete accord with Eq. (18). In the volume preserving case, we see that $\Lambda \Pi = 1$, and so the formulae Eq. (9), (12) and (14) again hold, and the average times are given by eq. (10) with $\lambda$ replaced by $\Lambda$.

Though the region $A$ that we considered above is special, a theorem of Easton [12] implies that the rate of escape for any isolating neighborhood of an invariant set are *asymptotically similar*. Thus any region surrounding the fixed point will have escape that is exponential with rate $\Lambda$. Similarly, any diagonalizable hyperbolic map with expansion $\Lambda$ will also have the same asymptotic rate.

We expect exponential decay of the transit time decomposition for a hyperbolic system. Numerical observations of transport, however, indicate that the transit time decomposition decays algebraically when there are elliptic regions in A. The simplest example of this behavior is the trivial shear:

$$x' = x + y'$$
$$y' = y$$

Let $A$ be the unit square as before. Now the entry set is the triangle $I = \{(x,y) = 0 \leq x < y \leq 1\}$. The transit time decomposition is

$$T_j = \left\{(x,y) \in I : 1 - jy < x < 1 - (j-1)y\right\} .$$

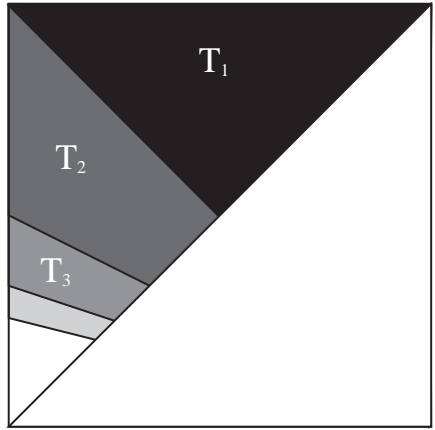

Figure 4. Transit time decomposition for the simple shear.

These sets have measures which decrease algebraically

$$\mu(T_1) = \frac{1}{4}, \quad \mu(T_j) = \frac{1}{j(j^2 - 1)} \quad j > 1$$

The sums to get the average transport times are elementary telescoping sums, and again these verify Eqs. (9) and (12). However, for this example the average transit time for $A$, (14), does not exist.

## 6 Bounded Orbits for the Hénon Map

As a final example we use the average transit time, which is straightfoward to compute, as an effective method to obtain the accessible area, which is not otherwise computable. Here we do this for the resonance zone of the Hénon map. The calculation involves several steps. First we find the points on the minimizing and minimax homoclinic orbits, $z_m$ and $z_h$, which bound the lobe, see Fig. 1. We then construct the boundary of the incoming set, by discretizing $W^u$ and $W^s$, as graphs $y^u(x)$ and $y^s(x)$, to a resolution $h = (x_m - x_h)/N$ for a fixed number of pixels $N$. Generally, in our calculation we used $N = 10^4$. By reversibility, the exit set is the reflection of the entrance set about $x+y=0$. The average exit time, $\langle t^+ \rangle_I$, is given by an integral of the piecewise constant exit time over the incoming set. We do this double integral in the most naive way by first integrating $T(x) = \int_{y^u(x)}^{y^s(x)} t^+(x,y)dy$ for a fixed $x$, and then integrating over $x$ using Simpson's rule. To compute $T(x)$, we use bisection to zoom in on the discontinuities of $t^+$: first evaluate $t^+$ on a grid of size $h$; if $t^+$ does not change between two grid points, we assume (possibly incorrectly) that it is constant between. If there is a change, we bisect the interval until either $t^+$ is equal on the endpoints, or the perceived error in neglecting the variation in $t^+$ is small enough (we chose an error of $10^{-3}$ for this). Then $T(x)$ is the sum of the $t^+$ values times the interval lengths. The resulting average exit time is shown in Fig. 5. While $\langle t^+ \rangle_I$ is generally decreasing, there are numerous small upward jumps. We discuss these more below.

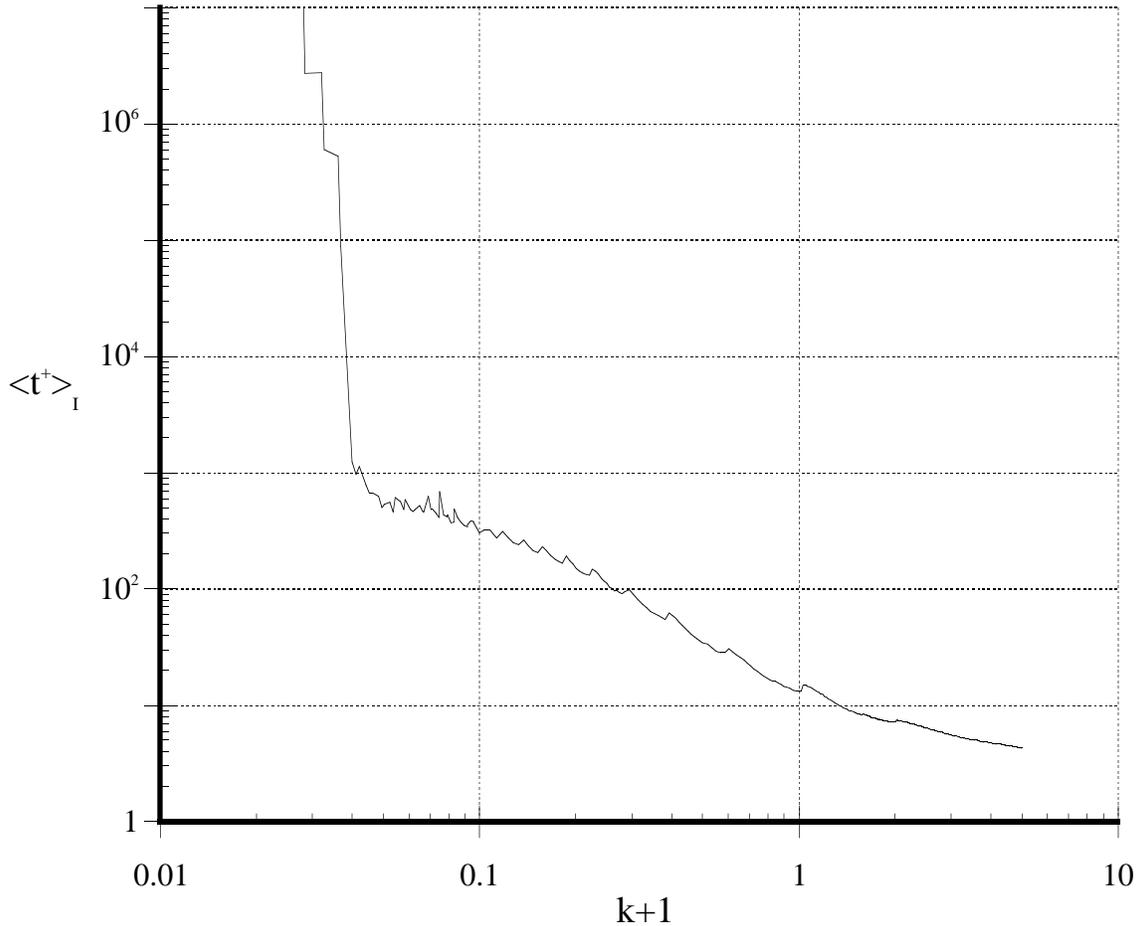

Figure 5. Average exit time for the Hénon map as a function of parameter. Most points used $N=10^4$ and $t_{max}=10^6$, though for k<-.95, we used $t_{max}=5\times10^7$

The area of the lobe is either given by summing the number of pixels in the exit set, or taking the difference in action between the two homoclinic points [24]. Then $\mu(A_{acc})$ is obtained from Eq.(11). It is interesting to compare this with the total area of the resonance zone itself, $\mu(A)$. This is most easily computed by taking the difference in action between the action of the minimax homoclinic point and the fixed point. These areas are shown in Fig. 6.

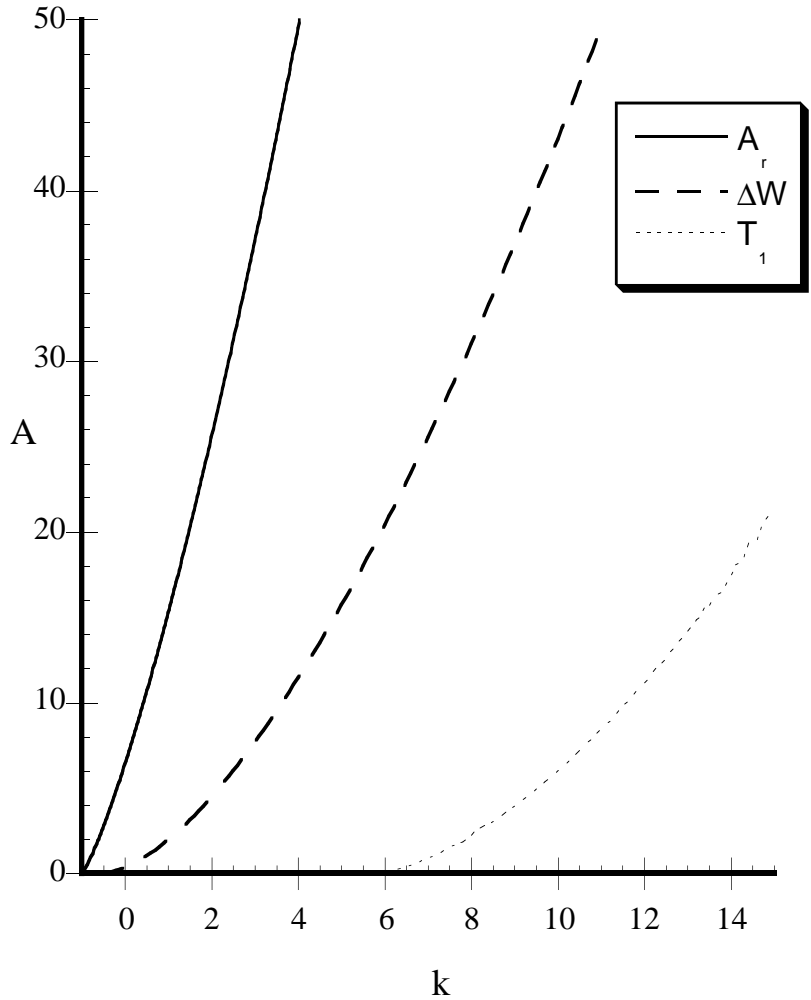

Figure 6. Resonance and Lobe area for the Hénon map. Also shown is $\mu(T_1)$, which is nonzero beyond the formation of the geometric horseshoe at $k \approx 5.706$.

Note that the resonance and lobe area grow monotonically and smoothly. The value $k = -1$ corresponds to the saddle node bifurcation, where the fixed point resonance zone is created. Slightly above this point, the lobe area, $\Delta W$ is apparently exponentially small and most of the resonance is filled with invariant curves.

Combining these results, using Eq. (11) gives the accessible area. In Fig. 7 we show the accessible fraction, $\mu(A_{acc})/\mu(A)$. The area that is inaccessible, which is identical to the measure of the bounded orbits is given from (11) by

$$\mu(A_i) = \mu(A) - \mu(A_{acc}) = \mu(A) - \mu(I)<t^+>_I \ ,$$

This area is shown in Fig 8, which is nearly identical to Fig 4. of [1], but our method allows us to compute the results to much higher accuracy. The accuracy can be seen better in the next figure, which is the best representation of this information, Fig. 9. This shows the inaccessible fraction, $\mu(A_i)/\mu(A)$.

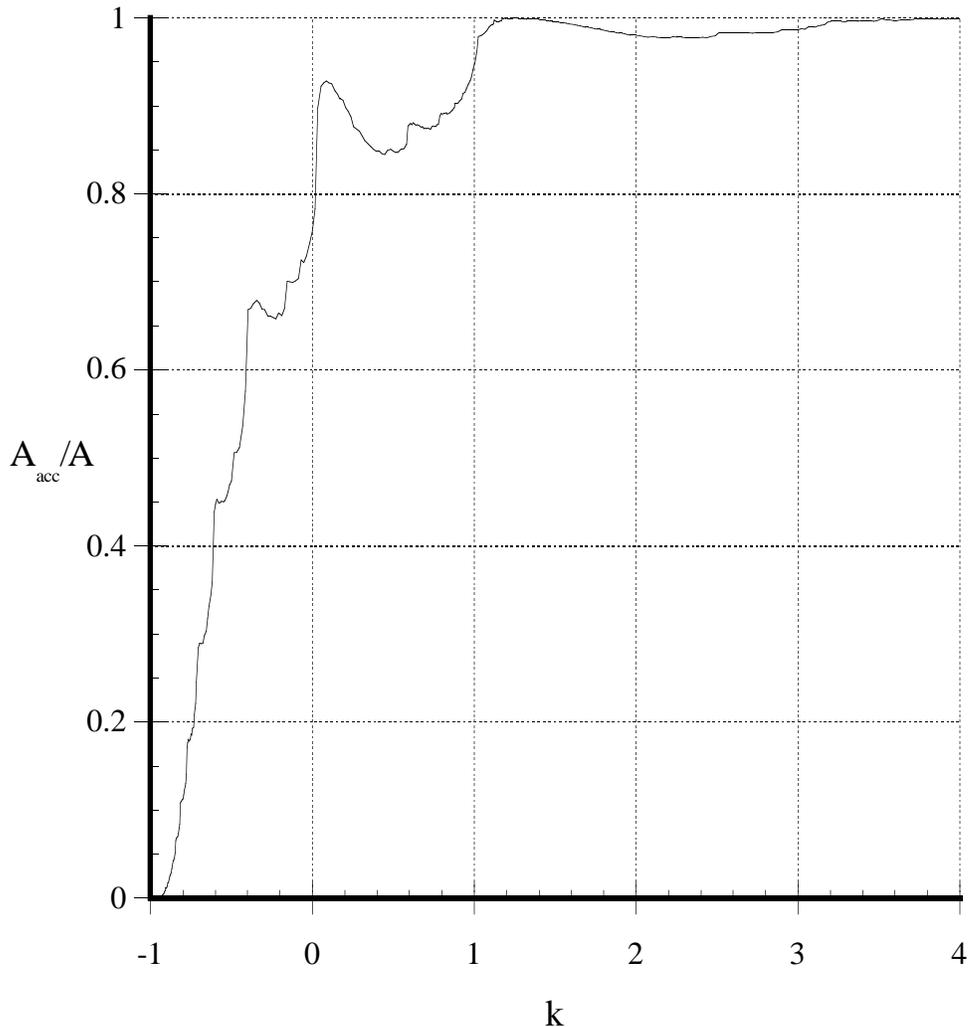

Figure 7. Accessible fraction for the fixed point resonance of the Hénon map.

In Fig. 9, the cut-off at a relative area of $10^{-3}$ is an artifact of our numerical method—it is quite difficult to reduce the error significantly. To do this one must increase the maximum number of iterations, $t_{max}$, to pick out narrow chaotic layers near invariant circles, and one must increase $N$ to find all possible discontinuities in the exit time. It would be nice improve the accuracy by detecting these discontinuities and find a form that fits $T(x)$ in their neighborhood, but we have not done that.

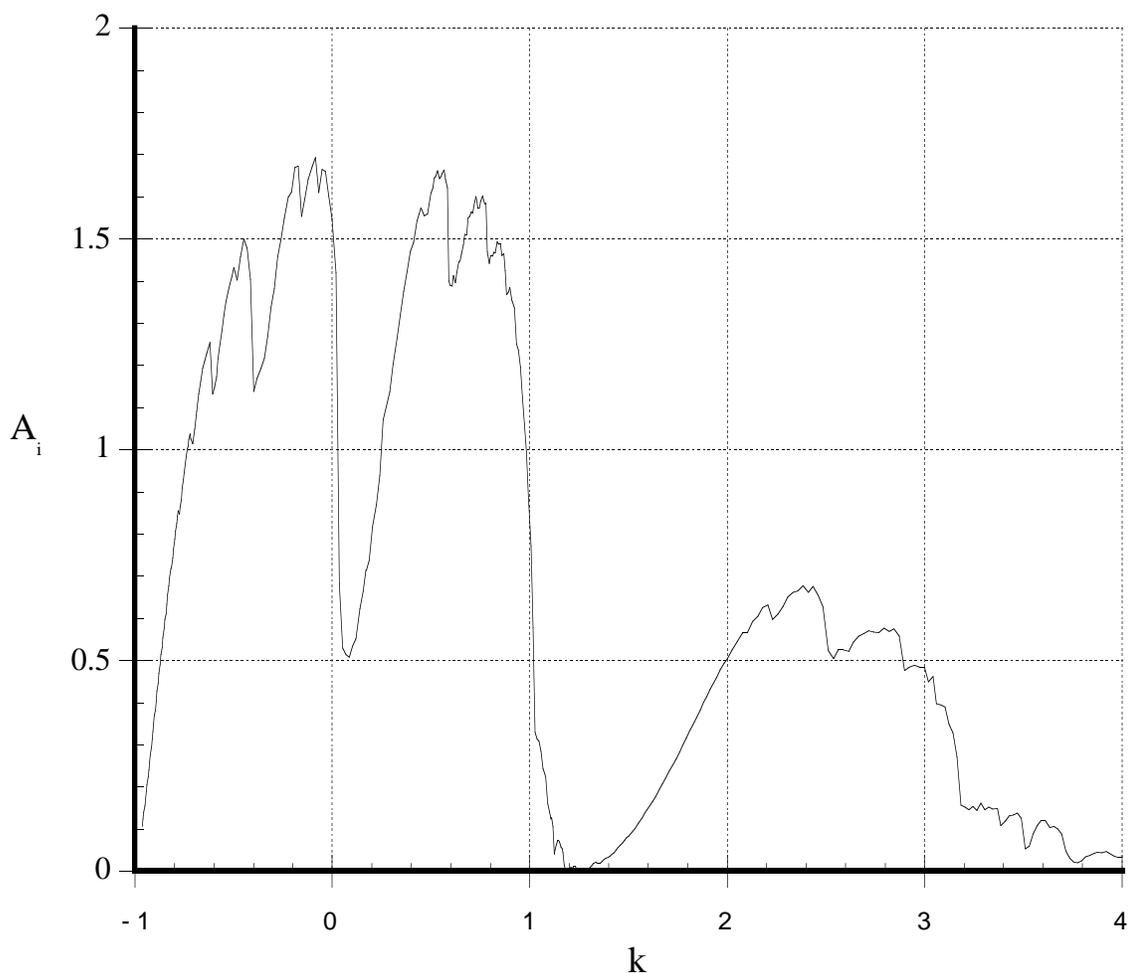

Fig. 8. Measure of the bounded orbits for the Hénon map.

The most prominent features in Fig. 9 are local minima near $k = 0$ and $k = 1.25$. These correspond to the quadrupling and tripling bifurcations of the elliptic fixed point. In the period 4 case, a pair of orbits with rotation number $1/4$ are created at the bifurcation. Interestingly for $0<k<0.4$, the saddle period 4 orbit, which has points on arranged on the square with corners ($\pm\sqrt{k}$, $\pm\sqrt{k}$), has nearly coincident stable and unstable manifolds—for practical purposes, a saddle connection, Fig. 10. Furthermore for $0.2<k<0.4$ this feature dominates the inaccessible set. Using this approximation gives $A_i \sim 4k$ for $k>0$—this approximation is shown in Fig 9, and gives good agreement with the computed results. Of course, there are other inaccessible islands, most importantly islands around the elliptic $\omega = 1/4$ orbit. These cause the fraction for $0<k<0.2$ to deviate from our simple form.

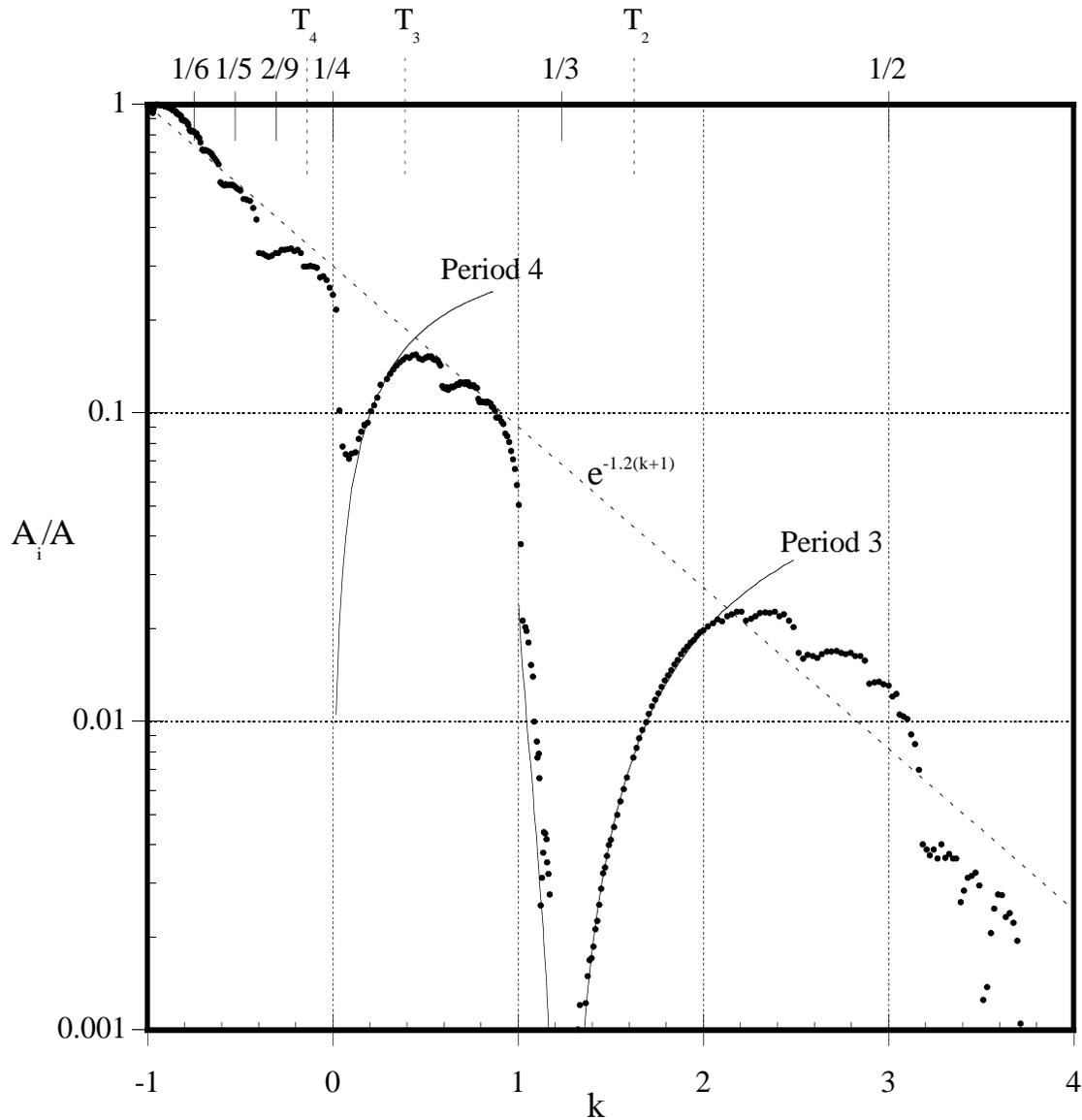

Fig 9. Inaccessible fraction of the resonance zone for the Hénon map. The two solid curves represent simple approximations to the area enclosed by the period 3 and period 4 saddle orbits. The fraction falls off on average exponentially with k, aside from dips near prominent bifurcation points. Along the top are shown bifurcations of the elliptic point (labeled by rotation number p/q), and homoclinic bifurcations corresponding to the creation of type 2, 3, and 4 trellises.

A pair of period three orbits is created by saddle node at $k=1$. At $k=5/4$ the period three saddle collides with the elliptic fixed point and there are no encircling invariant curves. Near this bifurcation, the most important feature is the virtually perfect saddle connection of the manifolds of the period three saddle (see Fig. 10), which has points at $(-\beta, \beta) \to (-1+\beta, \beta) \to (-\beta, 1-\beta)$, with $\beta = \sqrt{k-1}$. Using a triangle as an approximation for this area gives $A_i \sim \frac{1}{2}(2\beta-1)^2$. The resulting curve is also shown in Fig. 9; it fits remarkably well for $1.28 < k < 2$. Note that, contrary to the impression given by Fig. 9, the inaccessible fraction does not go to zero at $k=1.25$. In particular, there is a second period three orbit, $(\frac{1}{2}, -\frac{1}{2}) \to (-\frac{3}{2}, -\frac{1}{2}) \to (\frac{1}{2}, \frac{3}{2})$, that is elliptic. The island chain surrounding this orbit, see Fig. 11, has a relative area about $5(10)^{-4}$, accounting for most of the inaccessible area. There are no other islands visible.

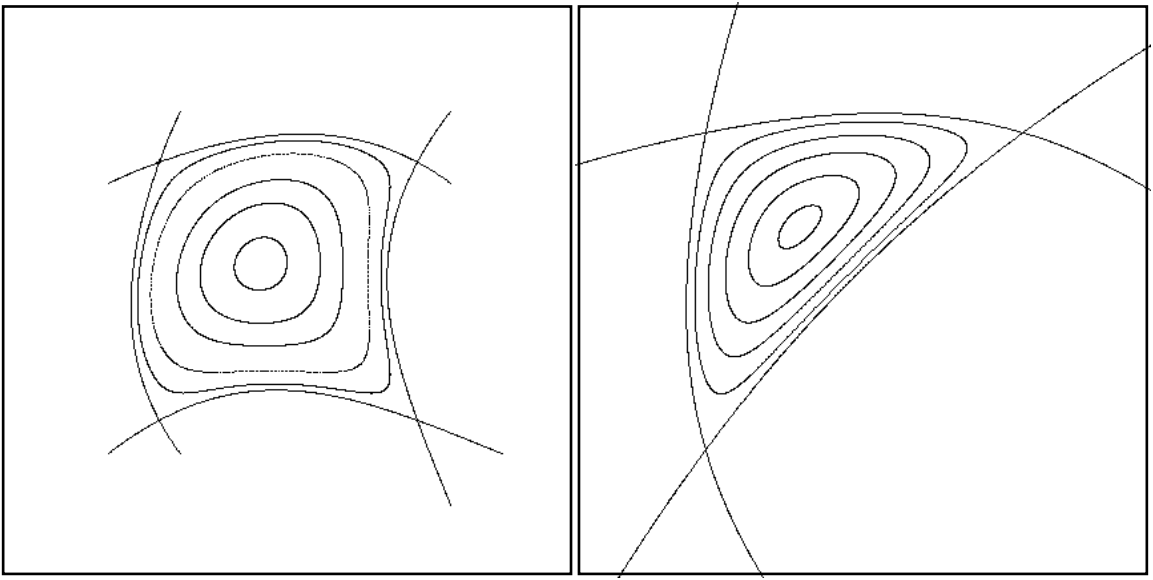

Fig. 10. Hénon map for $k = 0.2$ and $k = 1.6$. The outermost invariant curve surrounding the elliptic fixed point is closely approximated by a saddle connection of the $\omega = \frac{1}{4}$ and the $\frac{1}{3}$ saddle orbits, respectively. Bounds for the two figures are $(-1,1) \times (-1,1)$ and $(-1,0) \times (0,1)$, respectively

There are a number of other sharp drops in the inaccessible fraction, the most prominent occur at $k = -0.61, -0.413, 0.585, 0.78, 2.50$, and $3.17$. These occur when an invariant circle is destroyed, suddenly opening up a new accessible region. The newly opened region will be "large" if the critical circle is just outside a large island chain. For example for $k = -0.414$ there is an invariant circle just outside the $\omega = \frac{1}{5}$ island chain. This invariant circle is destroyed by $k = -0.413$, leading to the opening of a new accessible domain, and consequent decrease in measure of the bounded orbits.

## 7 Conclusions

We have shown that the average exit time from a region is given exactly by the ratio of the area of the accessible portion of the region to the area of the exit set. (11). It is interesting that this provides a justification for the oft used estimate that an "escape rate" from a region is given by the inverse of this ratio.

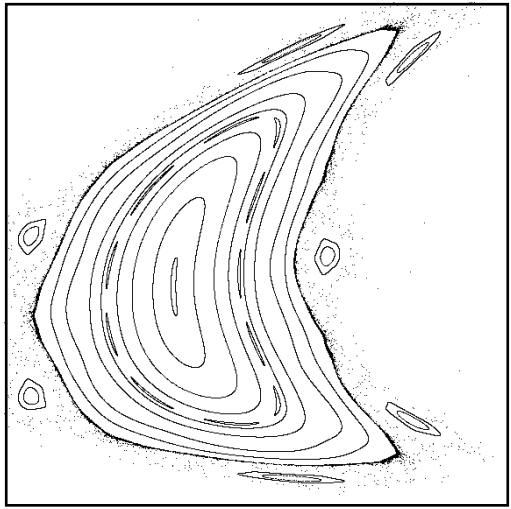

Figure 11 One of the islands around a period 3 orbit at *k* = 1.25, the period tripling point of the Hénon map. Plot bounds are (-1.51,-.597) × (-1.48,-.375).

We use this to provide a nice numerical tool for computing the measure of the bounded orbits for the Hénon map. Unfortunately, computational resources limit our accuracy in this calculation to a relative measure of about $10^{-3}$.

For the future, it would be nice to apply (11) to study the bounded orbits for higher dimensional maps, for example Moser's canonical form for the quadratic symplectic map [25]. Such maps are important in application to particle accelerators.

## Acknowledgments

I would also like to acknowledge the support of the National Science Foundation under grant #DMS-9305847, and support from a NATO Scientific Affairs Division Collaborative Research grant (#921181). This lead to discussions with Robert MacKay, whose insights were central to this research.